\documentclass{article} \usepackage{amssymb}  \newtheorem{lemma}{Lemma}
 \newcommand{\La}{\Lambda}
\newcommand{\R}{{\mathbb R}}  \newcommand{\Z}{{\mathbb Z}} 
 \newcommand{\C}{{\mathbb C}}  

\begin{document}
\begin{large}
\vspace{1cm}

\title{Approximation of discrete functions and size of spectrum\thanks{Published in Algebra i Analiz 21 (2009), no. 6, 227--240; translation in St. Petersburg Math. J. 21 (2010), no. 6, 1015--1025}}

\author{Alexander Olevskii\thanks{The first author is
 partially supported by the Israel Science Foundation} \ and Alexander Ulanovskii}

\date{} \maketitle

\noindent A.O.: School of Mathematics, Tel Aviv University\\ Ramat Aviv,  69978 Israel\\ E-mail:
 olevskii@post.tau.ac.il


\noindent A.U.: Stavanger University,  4036 Stavanger, Norway\\ E-mail: Alexander.Ulanovskii@uis.no




\begin{abstract}
   Let  $\La\subset\R$ be a uniformly discrete sequence and $S\subset\R$ a compact set.
   We prove that if there exists a bounded sequence of functions
   in Paley-Wiener space $PW_S$, which approximates $\delta-$functions
   on $\La$ with $l^2-$error $d$, then
measure$(S)\geq 2\pi (1-d^2) D^+(\La)$.
   This estimate is sharp for every $d$. Analogous estimate holds when the norms of approximating functions have a moderate growth, and we find a sharp growth restriction.
\end{abstract}

\medskip \medskip\noindent {\sl Keywords}: Paley--Wiener space; Bernstein space; Set of interpolation; Approximation of discrete functions

  \section{Introduction}

      \noindent {\bf 1.1.} Let $S$ be a compact set in $\R$, and let $m(E)$ denote the Lebesgue measure of $S$.
     By $PW_S$ we denote the Paley--Wiener space

$$
            PW_S: = \{f\in L^2(\R): f=\hat F, F=0 \mbox{ on } \R\setminus S\}
$$
      endowed with $L^2-$norm.
      Here $\hat F$ stands for the Fourier transform:
$$\hat F(x):=  \int_\R e^{itx}F(t)\,dt.$$

      By $B_S$ we denote the Bernstein space of bounded functions $f$
      (with the sup-norm), which are the Fourier transforms of Schwartz
      distributions supported by $S.$
     Clearly, every function $f\in PW_S$ (and every $f \in B_S$) can be extended to an entire function of finite exponential type.

     Throughout this paper we shall assume that $\La$ is a uniformly discrete set, that is
\begin{equation}
               \inf_{\lambda,\lambda'\in\La, \lambda\ne\lambda'} |\lambda-\lambda'|   >0.
\end{equation}
      The restriction operator
$$
                 f \to f|_\La
$$
      is a bounded linear operator from $PW_S$ into $l^2 (\La).$
      When this operator is surjective, the set $\La$ is called a set of
      interpolation   for $PW_S$.
      Similarly, if the restriction operator acts surjectively from $B_S$
      onto $l^\infty$, then $\La$ is called a set of interpolation  for $B_S$.
     The interpolation problem is to determine when $\La$ is a set of interpolation for $PW_S$ or $B_S$.

     The case $S= [a,b]$ is classical.
     Beurling and Kahane proved that in this case the answer  can be
     essentially given in terms of  the upper uniform density of $\La$,
         $$D^+(\Lambda):=\lim_{r\to\infty}\max_{a\in{\R}}\frac{\mbox{card}(\Lambda\cap(a,a+r))}{r}.
$$
  Namely, it was shown in \cite{k} that the condition
     $$m( S)> 2\pi     D^+(\La)
$$     is sufficient
    while the condition
\begin{equation}\label{1a}
    m(S)\geq 2\pi     D^+(\La)
\end{equation}
   is necessary for $\La$ to be a set of interpolation for $PW_S$.

    The first condition above is necessary and sufficient for $\La$ to be
    a set of interpolation  for $B_S$, see \cite{b1}.

     \medskip

     \noindent {\bf 1.2.} The situation becomes more delicate for the disconnected spectra.
   For the sufficiency part, not only the size but also the arithmetical structure
   of $\La$ is important.
   On the other hand,   Landau \cite{l} proved that (\ref{1a}) is still necessary
   for $\La$ to be a set of interpolation for $PW_S$, for every bounded set $S$.


   For compact spectrum $S$,  Landau's result can be
   stated in a more general form, which requires interpolation of $\delta-$functions only.
           For each $\xi \in \La$, let $\delta _\xi$ denote the corresponding $\delta$--function on $\La$:
$$
                     \delta_\xi (\lambda):=\left\{
                     \begin{array}{ll} 1 &  \lambda=\xi\\
                                   0 &  \lambda \neq \xi
                                   \end{array} \right. , \ \ \lambda\in\La.
$$

   \medskip

   \noindent{\bf Theorem A} (\cite{ou1}, Theorem 1)
     { \it Let $S$ be a compact. Suppose there exist functions $f_\xi\in PW_S$ satisfying
     $f_\xi|_\La=\delta_\xi, \xi\in\La,$ and \begin{equation}\label{norms}          \sup_{\xi\in\La} \Vert f_\xi\Vert < \infty.          \end{equation}
     Then inequality (\ref{1a}) holds.
     The statement is also true for $B_S-$spaces}.

\medskip

     \noindent
    {\bf 1.3.} The present paper is a direct continuation of \cite{ou1}.       We prove that the possibility of {\sl approximation} of $\delta-$functions on $\La$
       with a given $l^2 -$error already  implies an estimate from below on the measure of
       spectrum:

\medskip

\noindent{\bf
    Theorem 1} {\it Let $0<d<1,$ $S$ be a compact set, and $\La$ be a uniformly discrete set.
    Suppose there exist functions $f_\xi  \in PW_S$ satisfying (\ref{norms}) and such that

\begin{equation}\label{vertj}             \Vert f_\xi|_\La- \delta_\xi\Vert_{ l^2(\La)} \leq d,\  \mbox{ for every } \xi\in \La.     \end{equation}
      Then\begin{equation}\label{est}
                     m(S) \geq 2\pi (1-d^2) D^+(\La).
\end{equation}
    Inequality (\ref{est}) is sharp for every $d$.}

\medskip


     Clearly, by letting $d\to0$, Theorem 1 implies the necessary condition (\ref{1a}) for interpolation in $PW_S.$

            In sec. 4 we prove a version of Theorem 1 for the case when the norms of approximating functions
            have a moderate growth. On the other hand, no estimate of
            the measure of the spectrum    is possible if the norms grow too fast.

            In a contrast to Theorem 1 we show in sec. 5  that the possibility
            of $l^\infty$-approximation does not imply any estimate on the measure
            of $S$.     Similar problems for approximation in $l^p$ are discussed in sec. 6.

    Some  results of this paper were announced in \cite{ou}.

  \section{Lemmas}
Our approach to proof of Theorem 1 includes Landau's method (see \cite{l} and sec. 2 in \cite{ou1})  and some arguments from Kolmogorov's width theory.

\medskip\noindent {\bf 2.1. Concentration}

   \medskip\noindent {\bf  Definition}: Given a number $c ,0<c<1,$ we say that a
linear subspace $X$ of $L^2(\R)$ is
 {\it $c$-concentrated} on a set $Q$ if $$
                     \int_Q |f(x)|^2\, dx \geq c \Vert f\Vert_{L^2(\R)}^2, \ \mbox{ for every } \ f\in  X.
$$

\begin{lemma}\label{la} Given   sets $S,Q\subset\R$ of positive measure and a number $0<c<1$, let $X$ be a linear subspace of
$PW_{S}$ which is     $c$-concentrated on  $Q$.     Then
$$\mbox{dim}\, X \leq \frac{m(Q)\,m(S)}{2\pi c}.$$
\end{lemma}

\medskip
    This lemma is contained in \cite{l} (see statements
(iii) and (iv) in Lemma 1).

\medskip
\noindent
     {\bf 2.2. A remark on Kolmogorov's width estimate}

   \begin{lemma} Let $0<d<1$, and $\{{\bf u}_j\}, 1\leq j\leq n,$ be an orthonormal basis in an  $n$-dimensional complex
  Euclidean space $U$.
  Suppose that $\{{\bf v}_j\}, 1\leq j\leq n,$ is a family of vectors in $U$
  satisfying
\begin{equation}\label{10}
                \Vert {\bf v}_j-{\bf u}_j\Vert \leq d,\ j=1,...,n.
\end{equation}
  Then for every $\alpha$, $
              1<\alpha < 1/d,$
 there is a linear subspace $X$ in  $\C^n$ such that

         (i) dim$\, X >(1- \alpha ^2 d^2)n-1$;

         (ii) the estimate
$$ Q({\bf c}):=
\Vert \sum_{j=1}^n c_j {\bf v}_j\Vert^2\geq (1-\frac{1}{\alpha})^2\sum_{j=1}^n|c_j|^2,
$$
holds for every vector ${\bf c}=(c_1,...,c_n)\in X.$
\end{lemma}

\medskip
  \noindent
 The classical equality for Kolmogorov's width of "octahedron"
 (see \cite{ko}) implies that the dimension of the linear span of ${\bf v}_j$ is at least $(1-d^2)n$. This means that
  there exists a linear space $X\subset \C^n$, dim$\, X \geq (1-d^2)n$, such that
 the quadratic form $Q({\bf c})$
 is positive on the unite sphere of $X$.   Lemma 2 shows that by a
 small relative reduction  of the dimension, one can get an estimate of this
 form from below by a positive constant independent of $n.$

   We are indebted to E.Gluskin  for the following simple proof of this
 lemma.

\medskip\noindent {\bf Proof}.
   Given an $n\times n$ matrix $T=(t_{k,l}), k,l=1,...,n$, denote by $s_1(T)\geq...\geq s_n(T)$
   the singular values of this matrix (=the positive square roots of the eigenvalues of $TT^*$).

   The following properties are well--known (see  \cite{hj}, ch. 3):

\medskip
        (a)  (Hilbert--Schmidt norm of $T$ via singular values)   $$\sum_{j=1}^n s_j^2(T) = \sum_{k,l=1}^n |t_{ k,l}|^2.$$

\medskip
        (b) (Minimax--principle for singular values)  $$s_k(T) = \max_{L_k}  \min_{x\in L_k,\Vert x\Vert=1}\Vert Tx\Vert,$$where $\Vert\cdot\Vert$ is the norm in $\C^n$, and the maximum is taken over all linear subspaces $L_k\subseteq \C^n$ of dimension $k.$

   \medskip     (c)  $s_{k+j-1} (T_1+T_2) \leq s_k(T_1)+s_j(T_2)$, for all $ k,j\geq 1,k+j-1\leq n.$

\medskip
   Denote by $T_1$  the matrix, whose columns are the coordinates of ${\bf v}_l$ in
   the basis ${\bf u}_k$, and  set $T_2:=I-T_1$, where $I$ is the identity matrix.
     Then  property (a) and (\ref{10}) imply: $$\sum_{j=1}^n s_j^2(T_2) < d^2 n,$$
      and hence:$$
                   s^2_j(T_2) \leq  d^2\frac{ n}{j},\ 1\leq j\leq n.
$$
   This and (c) give:$$
           s_k(T_1) \geq s_n(I)-s_{n-k+1}(T_2)\geq 1- d\sqrt\frac{n}{n-k+1}  .
$$
 Since $s_n(I)=1$, by setting $k=n-[\alpha^2d^2n]$, where $[\cdot]$ means the integer part, we obtain:
 $$
           s_k(T_1) \geq 1-\frac{1}{\alpha}, \  k=n-[\alpha^2d^2n].
$$
 Now,  one can obtain  from (b)  that there exists $X$ satisfying the conclusions of the lemma.

\section{Proof of Theorem 1}
{\bf 3.1.}   Observe that  condition (\ref{norms}) implies
    the uniform boundedness of interpolating functions $f_\xi$:
\begin{equation}\label{1}
      |f_\xi (x)|  =\left| \int_{S} F_\xi(t)e^{ixt} dt\right| \leq \sqrt{m(S)} \Vert F_\xi\Vert_{L^2(\R)} < C_1.
\end{equation}

    We shall also  use  the following well--known fact (see \cite{ya}, Theorem 17): given a bounded spectrum $S$ and a uniformly discrete set $\La$,
    there exists $C(S,\La)$ such that
\begin{equation}\label{2}
  \sum_{\lambda \in \La} |f(\lambda)|^2 \leq C(S,\La) \int_{\R} |f(x)|^2\,dx, \mbox{ for every } f\in PW_S.
\end{equation}

\noindent
{\bf 3.2.}
Fix a small number $\delta>0$. Set $S(\delta):=S+[-\delta,\delta]$ and
\begin{equation}\label{varp}g_\xi(x):=f_\xi(x)\varphi(x-\xi),\ \xi\in\La, \
\varphi(x):=\left(\frac{\sin(\delta x/2)}{\delta x/2}\right)^2.
\end{equation}
 Clearly, $\varphi\in PW_{[-\delta,\delta]}$, so that $g_\xi\in PW_{S(\delta)}$. Also, since $\varphi(0)=0$ and $|\varphi(x)|\leq 1, x\in\R,$ it follows from (\ref{vertj})
that each $g_\xi|_\La$
approximates $\delta_\xi$ with an $l^2-$error $\leq d$:
\begin{equation}\label{gjl}
\Vert g_\xi|_\La-\delta_\xi\Vert_{l^2(\La)}\leq d, \ \xi\in\La.
\end{equation}

\medskip\noindent {\bf 3.3. }
 Fix numbers $a \in \R$ and $r>0$, and set
$$      I:=(a-r,a+r),\  \nu = \nu (I) := \mbox{card} (\La\cap I).$$
      From (1) we have:
\begin{equation}\label{3}
               \nu  < C |I|.
\end{equation}
   Here and below in this proof we denote by $C$ constants which do not depend on $I.$

    Denote by $\lambda_1<...<\lambda_\nu $ the elements of $\La \cap I$.
      It follows from (\ref{gjl}) that the vectors
$$
       {\bf v}_j:= (g_{\lambda_j} (\lambda_1),...,g_{\lambda_j}(\lambda_\nu))  \in \C^\nu, \ j=1,...,\nu,
$$
  satisfy (6) where $\{{\bf u}_j, j=1,...,\nu\}$ is the standard orthonormal
  basis in $\C^\nu$.

   Fix  a number $\alpha, 1<\alpha<1/d$.  By Lemma 2 there exists a subspace
$         X = X(a,r,\alpha)\subset \C^\nu$
    such that:

\medskip

            (i)    dim $X  > (1-\alpha^2d^2)\nu-1$,

\medskip
            (ii)   for every vector ${\bf c}=(c_1,c_2,...,c_\nu) \in X$
            the  inequality holds:
$$
       \Vert\sum_{j=1}^\nu c_j{\bf v}_j\Vert^2=    \sum_{k=1}^\nu\left|\sum_{j=1}^\nu c_j g_{\lambda_j} (\lambda_k)\right|^2
         \geq (1-\frac{1}{\alpha})^2\sum_{j=1}^\nu|c_j|^2.
$$
     Hence, we have from (\ref{2}) that
\begin{equation}\label{kk}
         \int_\R  \left|\sum_{j=1}^\nu c_jg_{\lambda_j}(x) \right|^2\,dx   \geq   C\sum_{j=1}^\nu |c_j|^2,   \  (c_1,...,c_\nu)\in X.
\end{equation}

\noindent {\bf 3.4.}   Set $I':= (a- r(1+\delta) , a+r(1+\delta)).$
       Then, due to (\ref{1}), (\ref{varp}) and (\ref{3}),  every function $$g(x):=\sum_{j=1}^\nu c_j g_{\lambda_j}(x)$$ satisfies:
$$
       \int_{\R\setminus I'}  |g(x)|^2\,dx  = \int_{\R\setminus I'}\left| \sum_{j=1}^\nu c_j f_{\lambda_j}(x)\left(\frac{ \sin \delta (x-\lambda_j)/2}{  \delta (x-\lambda_j)/2}\right)^2\right|^2\,dx$$$$
     \leq C \left(\sum_{j=1}^\nu |c_j|^2\right)   \int_{\R\setminus I'} \sum_{j=1}^\nu\frac{ 1}{\delta^4 (x-\lambda_j)^4}\,dx$$
     \begin{equation}\label{4}
     \leq C|I|\left(\sum_{j=1}^\nu |c_j|^2 \right) \frac{ 1}{\delta ^4} \int_{|y| >\delta r} \frac{ dy}{y^4}
     \leq \frac{C}{\delta^7r^2}  \sum_{j=1}^\nu |c_j|^2.
\end{equation}

      Fix $\epsilon>0$.   Inequalities (\ref{kk}) and (\ref{4}) show that there is a number $r_0=r(\delta, \epsilon)$
     (not depending on $a$ and ${\bf c}$)  such  that $r>r_0$ implies:
$$
              \int_{I'} |g(x)|^2\,dx  \geq (1-\epsilon) \int_{\R} |g(x)|^2\,dx.
$$
     This  means that the  subspace
$$
              G:= \{g(x)=\sum_{j=1}^\nu c_jg_{\lambda_j}(x); \ (c_1,...,c_\nu) \in X\}\subset L^2(\R)
$$
    is $(1-\epsilon)$-concentrated on $I'$, provided $r> r_0.$

\medskip

\noindent {\bf
3.5.}  Clearly, dim$\, G \geq$ dim$\, X$, so Lemma 1  now implies:$$
        \mbox{dim}\, X  \leq  \frac{m(S_\delta) |I'|}{2\pi(1-\epsilon)}  .$$
    Using  inequality (i) for dim$\,X$,  we obtain:
$$
       (1-\alpha^2 d^2)\nu -1 \leq  2r(1+\delta)\frac{m(S_\delta) }{2\pi(1-\epsilon) },
$$
  and so
      $$\frac{\mbox{card}\left( \La\cap (a-r,a+r)\right)}{2r}  \leq \frac{(1+\delta)m(S_\delta)}{2\pi (1-\epsilon)(1-\alpha^2d^2)}+\frac{1}{2r(1-\alpha^2d^2)}.$$

Now, for each fixed number $r$ we choose  $a$ so that  the left part is maximal,
    and then take limit as $r \to \infty$:
    $$
      D^+(\La)  \leq \frac{(1+\delta)\,m(S_\delta)}{2\pi(1-\epsilon)(1-\alpha^2d^2)}.
$$
     Since this inequality is true for all positive $\epsilon, \delta$ and every $\alpha\in (1,1/d) $,  we conclude that estimate (5) is true.

\medskip

\noindent {\bf
   3.6.}  Let us show  that  estimate (5) is sharp for every $d.$
  Pick up a number $a\in (0,\pi)$, and set $S:=[-a,a]$, $\La:=\Z$ and $$f_j(x):=\frac{\sin a(x-j)}{\pi (x-j)}\in PW_S, \ j\in\Z.$$ We have for every $j\in\Z$ that
$$
\Vert f_j|_\Z-\delta_j \Vert^2_{l^2(\Z)}=\Vert f_0|_\Z-\delta_0
\Vert^2_{l^2(\Z)}=\sum_{k\ne 0}\left(\frac{\sin a k}{\pi k}\right)^2
+\left(\frac{a}{\pi}-1\right)^2=
$$
$$
\frac{a}{\pi}-\frac{a^2}{\pi^2}
+\left(\frac{a}{\pi}-1\right)^2=1-\frac{a}{\pi}.
$$
Hence, the assumptions of Theorem 1 hold with
$d^2=1-a/\pi$. On the other hand, since $D^+(\Z)=1,$ we see that
$m(S)= 2\pi(1-d^2)D^+(\Z)$, so that estimate (\ref{est}) is
sharp.

\section{Moderate growth of norms}
{\bf 4.1.} Assume that the norms of functions in (\ref{vertj}) satisfy
\begin{equation}\label{restr}\Vert f_\xi\Vert_{L^2(\R)}\leq Ce^{|\xi|^\gamma}, \ \xi\in\La,\end{equation} where $C$ and  $\gamma$ are some positive constants. In this section we show that the statement of Theorem 1 remains true, provided  $\gamma<1$ and the density $D^+(\La)$ is replaced by the upper density $D^*(\La)$,
$$
D^{*}(\La):=\limsup_{a\to\infty}\frac{\mbox{card}\,(\Lambda\cap(-a,a))}{2a}.
$$
Restriction $\gamma <1$ is sharp: we show also that no estimate on the measure of spectrum is possible when the norms of $f_\xi$ grow exponentially.

Observe that  $D^*(\La)\leq D^+(\La)$, for each $\La$. However, one has $D^*(\La)= D^+(\La)$ whenever $\La$ is regularly distributed (in particular, when $\La$ is a bounded perturbation of integers).

\medskip

\noindent{\bf Theorem 2} {\it
 Let $0<d<1$.

 (i) Suppose $S$ is a compact set and  $\La$ is a uniformly discrete set. If there exist functions $f_\xi\in PW_S$ satisfying   (\ref{vertj}) and  (\ref{restr}) with some $0<\gamma<1$, then
 \begin{equation}\label{esti}
            m(S)\geq 2\pi(1-d^2)       D^*(\La). \end{equation}

(ii) For every $\epsilon>0$ there is a compact $S, m\,(S)<\epsilon$,
          a sequence $\La=\{n+o(1)\}$ and functions $f_\xi\in PW_S$ which satisfy
          (\ref{vertj}) and  (\ref{restr}) with $\gamma=1$. }


 \medskip\noindent{\bf Remark 1}.  Part (i) of Theorem 2 ceases to be true  for the density $D^+$, see  \cite{ou1}, Theorem 2.3.

 \medskip\noindent{\bf Remark 2}.  Similarly to \cite{ou1}, Theorem 2.4, one can check that the assumption $\gamma <1$ in part (i) can be weakened   by replacing
 it with any `non--quasianalytic' growth of norms in (\ref{restr}). It looks likely that the assumption $\gamma=1$ in part (ii) can be replaced with any `quasianalytic' growth. We leave this question open.

\medskip\noindent{\bf Remark 3}.  Let us show that if $S=[a,b]$ and $D^-(\La)\geq 1$, then assumption (\ref{vertj})  implies $b-a\geq 2\pi(1-d^2)$.
 Here $D^-$ is the lower uniform density of $\La$ (to define $D^-$, one replaces $\max$ with $\min$ in the definition of $D^+$).

 Recall that a set $\La$ is called a sampling set for $PW_S$
if there exist $A,B>0$ such that the inequality $$ A \Vert f\Vert_{L^2({\R})} \leq\left(\sum_{\lambda\in\Lambda}|f(\lambda)|^2\right)^{1/2} \leq  B \Vert
           f\Vert_{L^2({\R})}$$holds for every $f\in PW_S.$    The following   is a corollary of the classical result of Beurling on sampling sets in Bernstein spaces \cite{B}:
{\sl  Let $\La$ be a uniformly discrete set. If $D^-(\La)>a/\pi$ then  $\La$ is a sampling set for
$PW_{[-a,a]}$, if $D^-(\La)<a/\pi$, then it is not a sampling set for $PW_{[-a,a]}$.}

Now, suppose  $S=[a,b]$ and  $ D^-(\La)\geq 1$. Then $\La$ is a sampling set for  $PW_{[a,b]}$ provided $b-a<2\pi$. Clearly, in this case assumption (\ref{vertj}) implies  (3). By Theorem 1, we conclude that        $b-a \geq 2\pi (1-d^2).$

\medskip

Observe that $D^-(\La)=1$ for every $\La=\{n+o(1)\}$. It follows that the compact $S$ in part (ii) of Theorem 2 must be disconnected.
         On the other hand, we shall see that $S$ can be chosen a
        union of two intervals.


\medskip\noindent  {\bf 4.2. Proof of Theorem 2}

\medskip
 The  proof of part (i) is quite  similar to the proof of Theorem~1.

\medskip\noindent
{\bf 1.} Fix  numbers $\delta>0$ and $\beta,$ $\gamma<\beta<1$. There exists a function $\psi\in PW_{(-\delta,\delta)}$ with the properties:
\begin{equation}\label{ee}
\psi(0)=1,\ |\psi(x)|\leq 1, \ |\psi(x)|\leq C e^{-|x|^\beta},\ x\in\R,
\end{equation}
where $C>0$ is some constant.
It is well-known that such a function can be constructed as a product of $\sin(\delta_j x)/(\delta_jx)$ for a certain sequence  $\delta_j\to 0$.

Set$$h_\xi(x):=f_\xi(x)\psi(x-\xi), \ \xi\in\La.$$
Then each  $h_\xi$ belongs to $PW_{S(\delta)}$ and the restriction $h_\xi|_\La$
approximates $\delta_\xi$ with an $l^2-$error $\leq d$.

\medskip\noindent
{\bf 2.}   Set
$$
\Lambda_r:=\La\cap (-r,r),
$$
and denote by $C$ different positive constants independent on $r$.

 The argument in step 3.3  of the previous proof shows that there exists a linear space $X=X(r)$ of dimension $> (1-\alpha^2d^2)$card$(\La_r)-1$
such that
$$
\Vert \sum_{\xi\in\La_r}c_\xi h_\xi(x)\Vert^2_{L^2(\R)}\geq C\sum_{\xi\in\La_r}|c_\xi|^2,
$$
for every vector $(c_\xi)\in X$.

\medskip\noindent
{\bf 3.}   Since $\La$ is uniformly discrete, we have card$(\La_r)\leq Cr$. Further, using (\ref{restr}), similarly to (\ref{1}), we show that
$$
| f_\xi(x)|^2\leq Cm(S)\Vert  f_\xi\Vert_{L^2(\R)}^2\leq Ce^{C|\xi|^\gamma}\leq Ce^{Cr^\gamma}, \ \xi\in\La_r.
$$
These estimates  and (\ref{ee}) imply:
$$
\int_{|x|\geq r+\delta r}\left|\sum_{\xi\in\Lambda_r} c_\xi h_\xi(x)\right|^2\,dx=$$$$\int_{|x|\geq r+\delta r}\left|\sum_{\xi\in\Lambda_r}c_\xi
f_\xi(x)\psi(x-\xi)\right|^2\,dx \leq $$$$ \left(\sum_{\xi\in\Lambda_r}|c_\xi|^2\right)\left(Cr
e^{Cr^{\gamma}}\int_{|x|>\delta r}
e^{-2|x|^{\beta}}dx\right).
$$
Since $\beta>\gamma,$  the last factor  tends to zero as $r\to\infty.$ This and the estimate in step 4.2 show that  for every $\epsilon>0$ there exists $r_0=r(\delta,\epsilon) $ such that the linear space of functions
  $$\{h(x)=\sum_{\xi\in\La_r}c_\xi h_\xi(x); \ (c_\xi)\in X\}$$ is $(1-\epsilon)-$concentrated on $(-r-\delta r,r+\delta r)$, for all $r\geq r_0$. Moreover, the dimension of this space is at least $(1-\alpha^2d^2)\mbox{card}(\La_r)-1.$

\medskip\noindent
{\bf 4.} By Lemma 1, we obtain:
  $$
m(S(\delta))\geq  \frac{2\pi (1-\epsilon)}{1+\delta}  \frac{(1-\alpha^2 d^2)
(\mbox{card}(\La\bigcap (-r,r))-1)}{2r}.
   $$
Take now the upper limit as $r\to\infty$:     $$
m(S(\delta))\geq \frac{2\pi (1-\epsilon)}{1+\delta}(1-\alpha^2 d^2)D^*(\La).
   $$
Since this inequality holds for all $\epsilon>0,  \delta>0$ and $\alpha\in(1,1/d)$, we conclude that (\ref{esti}) is true.

\medskip\noindent {\bf 5.}  We shall now prove part (ii) of Theorem 2. We choose $S$ a union of two intervals and $\La$ a small perturbation of integers, as follows:
 $$S:=[-\pi-\epsilon,\pi+\epsilon]\cup[\pi-\epsilon,\pi+\epsilon], \ \La:=\{n+R^{-|n|-1}, n\in\Z\}.$$ Here $\epsilon>0$ is a given small number and $R>1$.

Denote by $\lambda_n:=n+R^{-|n|-1}$ the elements of $\La$, and set $$f_{\lambda_0}(x):=\frac{\sin \pi x}{\sin \pi\lambda_0}\cdot \frac{\sin\epsilon(x-\lambda_0)}{\epsilon(x-\lambda_0)},$$ and
$$
f_{\lambda_n}(x):=\frac{\sin \pi x}{\sin \pi\lambda_n}\cdot \frac{\sin \nu(n)(x-\lambda_n)}{\nu(n)(x-\lambda_n)}\cdot\prod_{|j|\leq 2|n|, j\ne n}\frac{\sin \nu(j)(x-\lambda_j)}{\sin \nu(j)(\lambda_n-\lambda_j)},n\ne0,
$$ where $\nu(n):=\epsilon/(4|n|+1)$. Observe that $m(S)=4\epsilon$, so to prove part (ii) it suffices to show that the functions $f_{\lambda_n}$ satisfy (4), provided $R$ is sufficiently large.

It is clear that $f_{\lambda_n}\in PW_S$, and that we have
\begin{equation}\label{tht}
 f_{\lambda_n}(\lambda_n)=1, \ n\in\Z, \ f_{\lambda_n}(\lambda_k)=0, |k|\leq 2n, k\ne n, n\ne 0.\end{equation} Further, we  assume that $R$ is  large enough so that the following three estimates hold for every $n\ne 0$ and every  $|k|>2|n|$:
$$
\left|\frac{\sin\pi\lambda_k}{\sin\pi\lambda_n}\right|\leq \frac{2\pi R^{-|k|-1}}{\pi R^{-|n|-1}}=2R^{-|k|+|n|};
$$
$$
\left|\frac{\sin \nu(n)(\lambda_k-\lambda_n)}{\nu(n)(\lambda_k-\lambda_n)}\right|\leq \frac{2}{\nu(n)(|k|-|n|)}\leq \frac{8}{\epsilon},
$$
and
$$
\left|\prod_{|j|\leq 2|n|, j\ne n}\frac{\sin \nu(j)(\lambda_k-\lambda_j)}{\sin \nu(j)(\lambda_n-\lambda_j)}\right|\leq\prod_{|j|\leq 2|n|, j\ne n}\frac{2}{ \nu(j)|j-n|}\leq
$$
$$
\left(\frac{2}{\nu(2n)}\right)^{4|n|}\frac{1}{ |n|!(3|n|)!}\leq \left(\frac{C}{\epsilon}\right)^{4|n|},
$$
 where $C>1$ is an
  absolute constant.
These estimates yield:
$$
|f_{\lambda_n}(\lambda_k)|\leq
16\left(\frac{C}{\epsilon}\right)^{4|n|+1}R^{-|k|+|n|}, \ |k|>2|n|, \ n\ne0.
$$ A similar estimate holds for $f_{\lambda_0}(\lambda_k)$ for each $k\ne 0$. Clearly, these estimates and (\ref{tht}) prove (4), provided $R$  is large enough.

  \section{$l^\infty-$approximation.}

  {\bf 5.1.} In a sharp contrast to Theorem 1, the possibility of $l^\infty$--approximation of $\delta$--functions on $\La$
   does not imply any restrictions on  the measure of spectrum.

    For approximation by $PW-$functions this follows from Lemma~3.1 in \cite{ou1}:
{\sl    For every  $N\geq 2$ there exists a set $S(N)\subset(-N,N)$, $m(S(N))=\frac{2}{N},$ such that
   $$
   \left|\frac{N}{2}\int_{S(N)}e^{itx}\,dt-\frac{\sin Nx}{Nx}\right|\leq \frac{C}{N}, \ x\in\R,
   $$
   where $C>0$ is an absolute constant.}

   The function $\sin Nx/Nx$ is essentially localized in a small neighborhood of the origin, and its Fourier transform is
   the unite mass uniformly distributed over the interval $[-N,N]$. The lemma shows that one can
re-distribute this mass over a set of small measure so that the `uniform error' in the Fourier transform is $O(1/N)$.

     For the $B_S-$functions, the result can be stated even in a stronger
     form:

\medskip

\noindent{\bf
    Proposition 1}  {\it Given a number $0<d<1$ and a uniformly discrete set $\La$, there exist a compact set
    $ S$ of measure zero and a bounded sequence of functions $f_\xi \in B_S$ satisfying
$$                       \Vert f_\xi|_\La - \delta_\xi \Vert_{l^\infty(\La)} \leq d,   \ \mbox{ for every } \xi \in \La.$$
     The set S can be chosen depending only on $d$ and the separation constant  in (1).}

\medskip


      Let us invoke the classical Menshov example from the
       uniqueness theory            of trigonometric series. It can be stated as follows (see \cite{bari} ch.14, sec.12, and remark in
      sec.18): {\sl
     There is a singular probability measure $\nu$ with compact support,
     such that
$$
                      \hat \nu(x) \to 0,\ \   |x|\to \infty.
$$
}

\medskip\noindent{\bf    Corollary 1} {\it  For every $\epsilon>0$  there is a compact set $S\subset\R$ of Lebesgue
    measure zero     and a function $f \in B_S$, such that
$$
                  f(0) = \Vert f\Vert_{L^\infty(\R)}  =1,\ \mbox{ and }  \                  |f(t)| < \epsilon,  \ \  |t|>\epsilon.
$$}
Indeed, it suffices to  set $f(x)=\hat\nu(cx),$  where $c$ is sufficiently large.

     Now Proposition 1 follows immediately:         take a positive number $\epsilon<\min\{d, \gamma(\La)\}$, where  $\gamma(\La)$ is defined in (1). Let $f$ be a function from the colollary. Then  the functions  $f_\xi(x):=f(x-\xi), \xi\in\La,$ satisfy the assumptions of Proposition~1.

\medskip

\noindent {\bf 5.2.}  Notice that the Bernstein space $B_S$  can be defined
    in a similar way  for every unbounded closed spectrum $S$ of finite measure, see \cite{ou1}.
    In \cite{ou1} we  constructed unbounded spectra $S$ of
    arbitrarily small measure  such that every uniformly discrete set $\La$ is a set of interpolation for $B_S$.
     This  was done by a certain iteration argument, using Lemma 3.1
     from that paper.    Using instead Corollary 1, one can prove by the same approach a more precise version of the result:

    \medskip

    \noindent{\bf Theorem 3}
{\it
     There is a closed set $S$ of measure zero such that every
   uniformly discrete set $\La$ is a set of interpolation  for $B_S$.
}

 \medskip
     \noindent
 {\bf  Remark}      Assumption $m(S)=0$ in Proposition 1 and Theorem 3 can be replaced
     by a stronger metrical `thinness' condition:  $S$ may have measure zero
      with respect to any given Hausdorff
     scaling function. 
For such an improvement one needs
     to use measures $\nu$ constructed in \cite{im}.

  \section{$B_S^p-$spaces and  $l^p-$approximation }

    One can include spaces $PW_S$ and $B_S$
    into a continuous chain of Banach spaces:
    Given a compact set $S$  and a number $p ,1\leq p\leq\infty$,
    denote by  $B_S^p$ the space of all entire functions $f\in L^p(\R)$
    that can be represented as the Fourier
    transform of a distribution $F$ supported by $S$.
        Clearly,  $B_S^2=PW_S$ and $ B_S^\infty= B_S$.

     Observe that for $p<p'$, one has the embedding
     \begin{equation}\label{emm}
              B^p_S\subset  B^{p'}_S
              \end{equation}
     with the corresponding inequality for norms.

     Let $\La$ be a uniformly discrete set. It is well-known that the restriction operator
        $f\to f|_\La$ acts boundedly from $B^p_S$ into $l^p(\La)$ (see, for example, \cite{ya}, p.82).
      $\La$ is called a set of interpolation  for $B_S^p$
     if this operator is surjective.

     Theorem A implies:

     \medskip\noindent{\bf Theorem 4}
    {\it Let $S$ be a compact and $p\geq 1$. If there exist functions $f_\xi \in B^p_S$ satisfying $f_\xi|_\La=\delta_\xi,\xi\in\La,$ and
    \begin{equation}\label{emn}
           sup_{\xi\in\La} \Vert f_\xi\Vert_{L^p(\R)}<\infty,
      \end{equation}
    then  condition (2) holds}.

\medskip
     In particular, this shows that if $\La$ is a set of  interpolation for $B^p_S$,
     then estimate (2) is true.

     However, when considering   $l^p(\La)-$approximation
     by functions from $B^p_S$,
     one should distinguish between the following two cases:
             $1\leq p\leq2$  and  $  2<p\leq\infty .$
     In the first case,  the measure of spectrum
     admits an estimate from below as in Theorem 1, while in the second case
     it does not as in Proposition 1:

\medskip

\noindent{\bf Theorem 5} {\it Let $0<d<1$  and $\La$ be a uniformly discrete set.

\medskip
(i) Suppose  $1\leq p\leq 2$ and $S$ is a compact set. If  every $\delta_\xi, \xi\in\La,$ admits approximation \begin {equation}\label{enn}\Vert f_\xi|_\La-\delta_\xi\Vert_{l^p(\La)} \leq d,\ \xi\in\La, \end{equation} by  functions $f_\xi\in B^p_S$ satisfying (\ref{emn}),   then condition (5)  holds true.

\medskip
(ii) Suppose $ p>2$. There exist a compact set $S\subset\R$ of measure zero and  functions $f_\xi\in B^p_S$ satisfying (\ref{emn}) and (\ref{enn}).
}

\medskip

 Part (i) is a consequence of Theorem 1, embedding (\ref{emm}) and the standard inequality
       between $l^p$ norms.

       Part (ii) follows  form the refinement  of Menshov's example (see \cite{ivm}):
      {\it There is a singular measure $\nu$ with compact support satisfying }
       $$
                   \hat\nu(x) = O( |x|^{-1/2}), \ \ |x|\to\infty.
$$

\medskip\noindent{\bf
Acknowledgment}: A part of the present work was done at the Mathematisches
Forschungsinstitut Oberwolfach during a stay within
the Research in Pairs Programme, April 2009. The authors  appreciate the hospitality  of the Institute.

\end{large}
\end{document}